\documentclass[letterpaper]{article}
\usepackage[margin=1in,dvips]{geometry}
\usepackage{graphicx,psfrag,amsmath,amsthm,amssymb}
\usepackage{natbib}
\usepackage{url}

\newcommand{\A}{\ensuremath{\mathbf{A}}}

\newcommand{\C}{\ensuremath{\mathbf{C}}}

\newcommand{\I}{\ensuremath{\mathbf{I}}}

\newcommand{\PP}{\ensuremath{\mathbf{P}}}

\newcommand{\RR}{\ensuremath{\mathbf{R}}}

\renewcommand{\b}{\ensuremath{\mathbf{b}}}

\newcommand{\dd}{\ensuremath{\mathbf{d}}}

\newcommand{\bl}{\ensuremath{\mathbf{l}}}

\newcommand{\q}{\ensuremath{\mathbf{q}}}

\newcommand{\uu}{\ensuremath{\mathbf{u}}}
\newcommand{\vv}{\ensuremath{\mathbf{v}}}

\newcommand{\x}{\ensuremath{\mathbf{x}}}
\newcommand{\y}{\ensuremath{\mathbf{y}}}
\newcommand{\z}{\ensuremath{\mathbf{z}}}
\newcommand{\0}{\ensuremath{\mathbf{0}}}

\newcommand{\blambda}{\ensuremath{\boldsymbol{\lambda}}}

\newcommand{\bzeta}{\ensuremath{\boldsymbol{\zeta}}}

\newcommand{\bbR}{\ensuremath{\mathbb{R}}}

\newcommand{\calL}{\ensuremath{\mathcal{L}}}

\newcommand{\calO}{\ensuremath{\mathcal{O}}}

\newcommand{\norm}[1]{\left\lVert#1\right\rVert}

{%
\begin{list}{#1}{
\vspace{-\topsep}
\vspace{-\partopsep}
\setlength{\itemindent}{0cm}
\setlength{\rightmargin}{0cm}
\setlength{\listparindent}{0cm}
\settowidth{\labelwidth}{#1}
\setlength{\leftmargin}{\labelwidth}
\addtolength{\leftmargin}{\labelsep}
\setlength{\itemsep}{0cm}
}%
}%
{%
\end{list}
\vspace{-\topsep}
\vspace{-\partopsep}
}

{\begin{enumerate}%
}%
{\end{enumerate}}

\newcommand{\Domain}{\operatorname{dom}}
\newcommand{\dom}[1]{\ensuremath{\Domain\left(#1\right)}}

\newcommand{\proxop}{\operatorname{prox}}
\newcommand{\prox}[2]{\ensuremath{\proxop_{#1}\left(#2\right)}}

\theoremstyle{plain}

\newtheorem*{lemma*}{Lemma}

\newtheorem*{prop*}{Proposition}

\theoremstyle{definition}

\newtheorem*{defn*}{Definition}

\newtheorem*{exmp*}{Example}

\newtheorem*{conj*}{Conjecture}

\theoremstyle{remark}

\newtheorem*{rmk*}{Remark}

\graphicspath{{grf/}}

\bibpunct[, ]{(}{)}{;}{a}{,}{,} 

\title{An ADMM algorithm for solving \\ a proximal bound-constrained quadratic program}
\author{
  Miguel \'A. Carreira-Perpi\~n\'an \\
  Electrical Engineering and Computer Science, University of California, Merced \\
  {\url{http://eecs.ucmerced.edu}}
}
\date{December 29, 2014}

\begin{document}

\maketitle

\begin{abstract}

  We consider a proximal operator given by a quadratic function subject to bound constraints and give an optimization algorithm using the alternating direction method of multipliers (ADMM). The algorithm is particularly efficient to solve a collection of proximal operators that share the same quadratic form, or if the quadratic program is the relaxation of a binary quadratic problem.

\end{abstract}

\section{Introduction}

We consider the convex \emph{proximal bound-constrained quadratic program (QP)}:
\begin{equation}
  \label{e:proxbqp}
  \min_{\x\in\bbR^D}{ f(\x) + \frac{\mu}{2} \norm{\x - \vv}^2 } \text{ s.t.\ } \bl \le \x \le \uu \qquad f(\x) = \frac{1}{2} \x^T \A \x - \b^T \x
\end{equation}
where $\mu > 0$, $\vv,\bl,\uu,\b\in\bbR^D$ and the $D \times D$ matrix \A\ is symmetric positive definite or semidefinite. The problem has the form of a convex proximal operator \citep{CombetPesquet11a,Moreau62a,Rockaf76b} $\x = \prox{f/\mu}{\vv}$ (with $\dom{f} = [\bl,\uu]$), and a unique solution. We will also be interested in solving a collection of such QPs that all have the same matrix \A\ but different vectors \b, \vv, \bl\ and \uu.

It is possible to solve problem~\eqref{e:proxbqp} in different ways, but we want to take advantage of the structure of our problem, namely the existence of the strongly convex $\mu$ term, and the fact that the $N$ problems have the same matrix \A. We describe here one algorithm that is very simple, has guaranteed convergence without line searches, and takes advantage of the structure of the problem. It is based on the alternating direction method of multipliers (ADMM), combined with a direct linear solver and caching the Cholesky factorization of \A.

\paragraph{Motivation}

Problem~\eqref{e:proxbqp} arises within a step in the binary hashing algorithm of \citet{CarreirRaziper15a}. The step involves $N$ independent subproblems
\begin{equation*}
  \min_{\x_n\in\{0,1\}^D}{ \frac{1}{2} \norm{\C\x_n - \dd_n}^2 + \frac{\mu}{2} \norm{\x_n - \vv_n}^2 } \qquad n = 1,\dots,N.
\end{equation*}
Intuitively, subproblem $n$ tries to map linearly a binary vector $\vv_n$ onto a real vector $\dd_n$ with minimal error, given a mapping of matrix \C. Hence, each subproblem is of the form of~\eqref{e:proxbqp}, but where $f$ is a least-squares function (where the rectangular matrix \C\ is common to all subproblems) and the variables are binary. Relaxing these subproblems results in $N$ proximal bound-constrained QPs of the form of~\eqref{e:proxbqp}.

\section{Solving a QP using ADMM}

We briefly review how to solve a quadratic program (QP) using the alternating direction method of multipliers (ADMM), following \citet{Boyd_11a}. Consider the QP
\begin{align}
  \label{e:qp}
  \min_{\x} & \quad \frac{1}{2} \x^T \PP \x + \q^T\x \\
  \text{s.t.} & \quad \A\x = \b,\ \x \ge \0
\end{align}
over $\x \in \bbR^D$, where \PP\ is positive (semi)definite. To use ADMM, we first introduce new variables $\z \in \bbR^D$ so that we replace the inequalities with an indicator function $g(\z)$, which is zero in the nonnegative orthant $\z \ge \0$ and $\infty$ otherwise. Then we write the problem as
\begin{align}
  \label{e:qp-admm}
  \min_{\x} & \quad f(\x) + g(\z) \\
  \text{s.t.} & \quad \x = \z
\end{align}
where
\begin{equation*}
  f(\x) = \frac{1}{2} \x^T \PP \x + \q^T\x,\qquad \dom{f} = \{\x\in\bbR^D\mathpunct{:}\ \A\x = \b\}
\end{equation*}
is the original objective with its domain restricted to the equality constraint. The augmented Lagrangian is
\begin{equation}
  \label{e:auglag}
  \calL(\x,\z,\y;\rho) = f(\x) + g(\z) + \y^T (\x-\z) + \frac{\rho}{2} \norm{\x-\z}^2
\end{equation}
and the ADMM iteration has the form:
\begin{align*}
  \x &\leftarrow \arg\min_{\x}{\calL(\x,\z,\y;\rho)} \\
  \z &\leftarrow \arg\min_{\z}{\calL(\x,\z,\y;\rho)} \\
  \y &\leftarrow \y + \rho(\x - \z)
\end{align*}
where \y\ is the dual variable (the Lagrange multiplier estimates for the constraint $\x = \z$), and the updates are applied in order and modify the variables immediately. Here, we use the scaled form of the ADMM iteration, which is simpler. It is obtained by combining the linear and quadratic terms in $\x-\z$ and using a scaled dual variable $\bzeta = \y/\rho$:
\begin{align*}
  \x &\leftarrow \arg\min_{\x}{\left( f(\x) + \frac{\rho}{2} \norm{\x - \z + \bzeta}^2 \right)} \\
  \z &\leftarrow \arg\min_{\z}{\left( g(\z) + \frac{\rho}{2} \norm{\x - \z + \bzeta}^2 \right)} \\
  \bzeta &\leftarrow \bzeta + \x - \z.
\end{align*}
Since $g(\x)$ is the indicator function for the nonnegative orthant, the solution of the \z-update is simply to threshold each entry in $\x+\bzeta$ by taking is nonnegative part. Finally, the ADMM iteration is:
\begin{subequations}
  \label{e:admm}
  \begin{align}
    \label{e:admm1}
    \x &\leftarrow \arg\min_{\x}{\left( f(\x) + \frac{\rho}{2} \norm{\x - \z + \bzeta}^2 \right)} \\
    \label{e:admm2}
    \z &\leftarrow (\x + \bzeta)_+ \\
    \label{e:admm3}
    \bzeta &\leftarrow \bzeta + \x - \z
  \end{align}
\end{subequations}
where the updates are applied in order and modify the variables immediately, $(t)_+ = \max(t,0)$ applies elementwise, and $\norm{\cdot}$ is the Euclidean norm. The penalty parameter $\rho > 0$ is set by the user, and $\z = (z_1,\dots,z_D)^T$ are the Lagrange multiplier estimates for the inequalities. The \x-update is a quadratic objective whose solution is given by a linear system. The ADMM iteration consists of very simple updates to the relevant variables, but its success crucially relies in being able to solve the \x-update efficiently. Convergence of the ADMM iteration~\eqref{e:admm} to the global minimum of problem~\eqref{e:qp} in value and to a feasible point is guaranteed for any $\rho > 0$.

\section{Application to our QP}

\begin{figure}[t]
  \centering
  \psfrag{z}[][]{$z$}
  \psfrag{Mz}[][]{$M(l,u,z)$}
  \psfrag{0}[B][B]{\raisebox{-5pt}{$l$}}
  \psfrag{1}[B][B]{\raisebox{-5pt}{$u$}}
  \psfrag{00}[B][B]{$l$}
  \psfrag{01}[B][B]{$u$}
  \begin{tabular}{@{}c@{\hspace{0.01\linewidth}}c@{\hspace{0.01\linewidth}}c@{\hspace{0.03\linewidth}}c@{}}
    \includegraphics[width=0.23\linewidth]{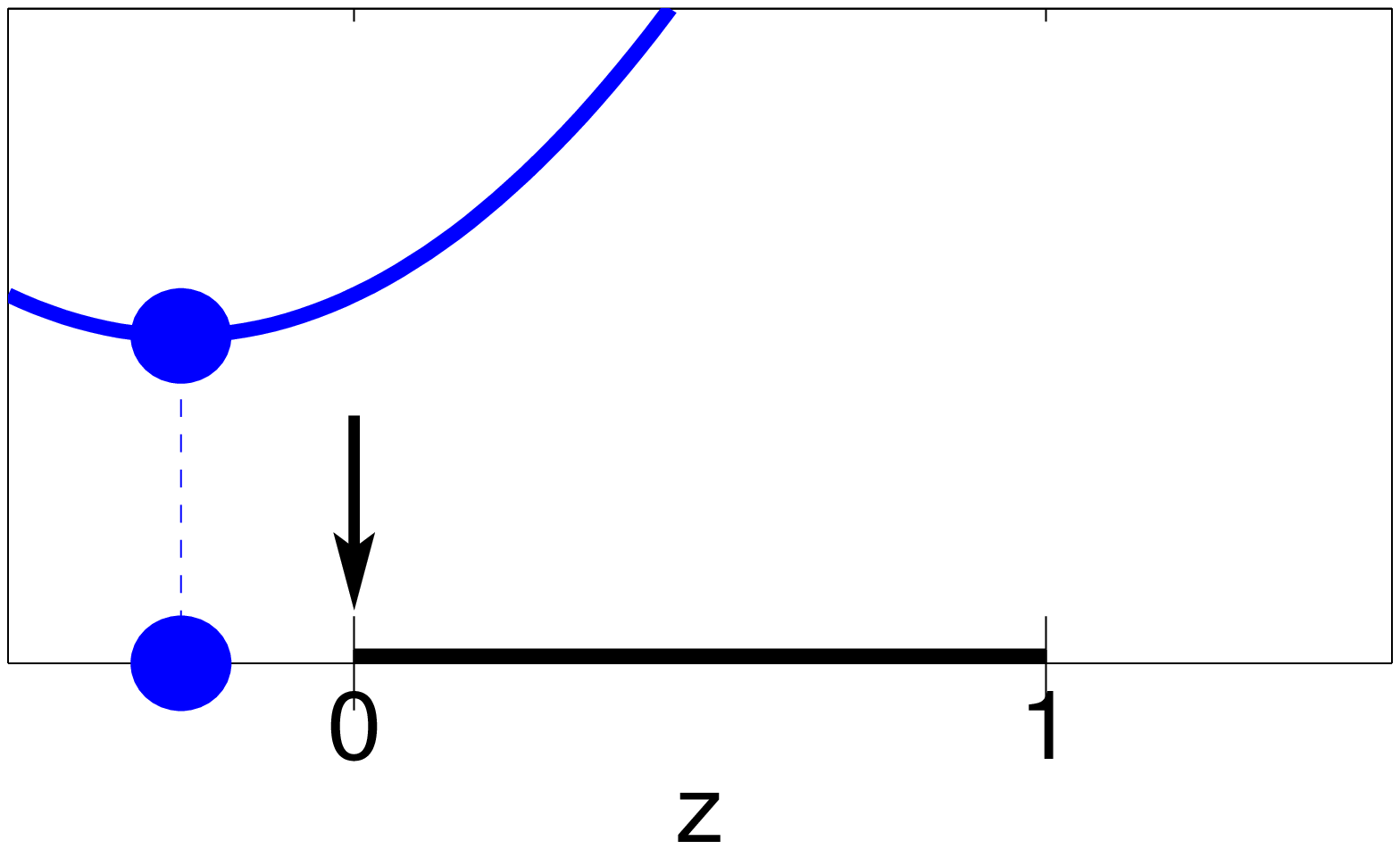} &
    \includegraphics[width=0.23\linewidth]{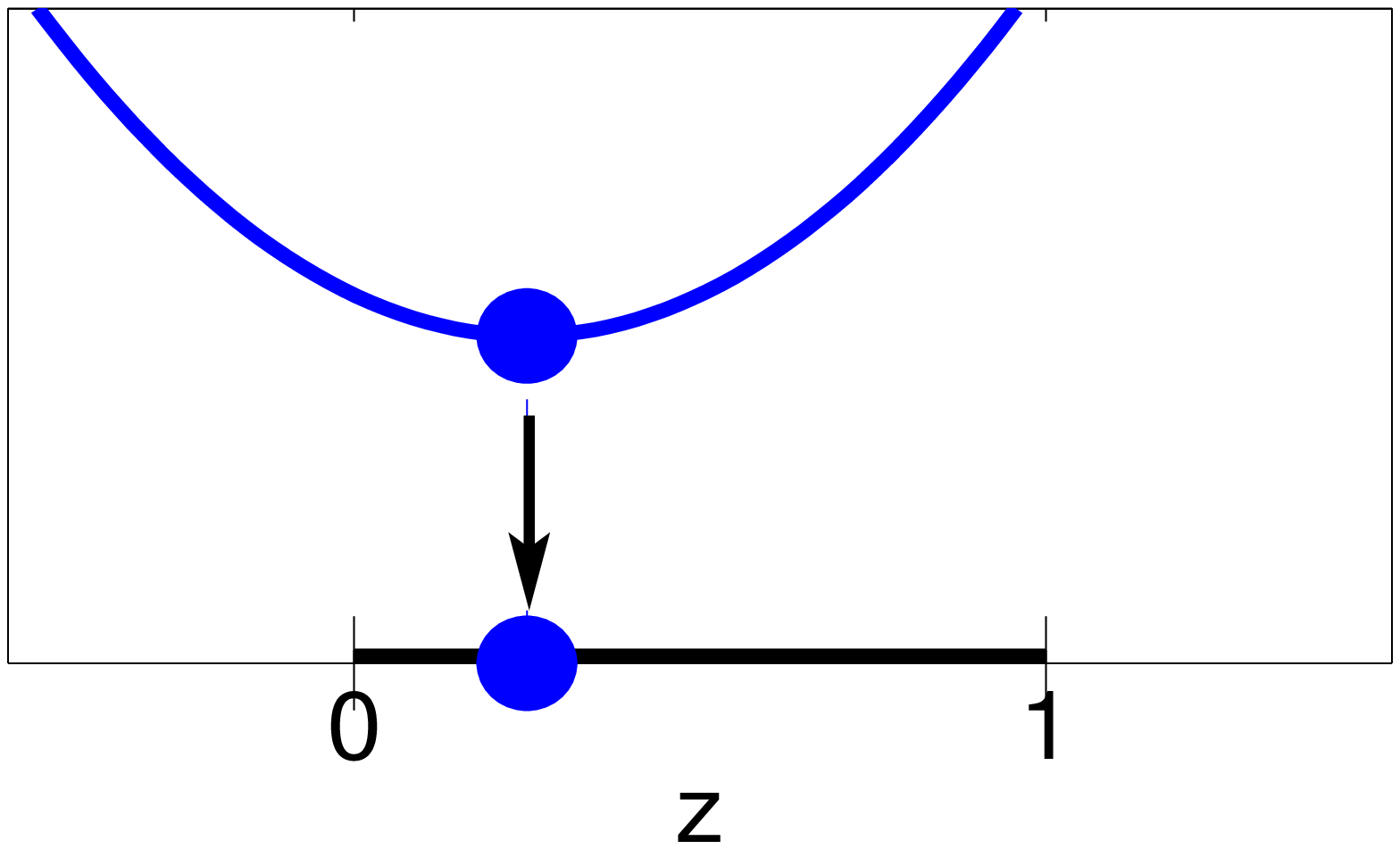} &
    \includegraphics[width=0.23\linewidth]{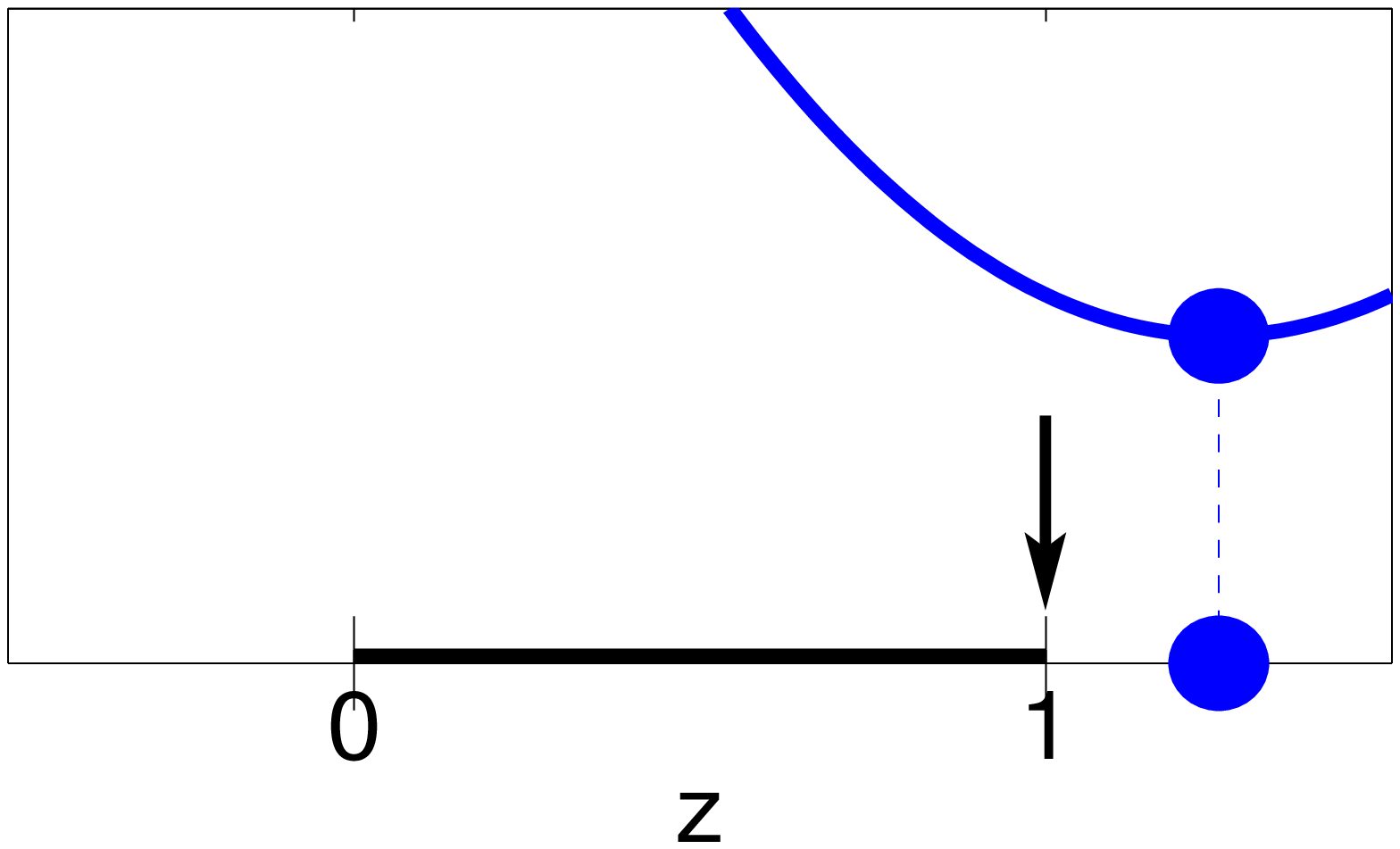} &
    \includegraphics[width=0.25\linewidth]{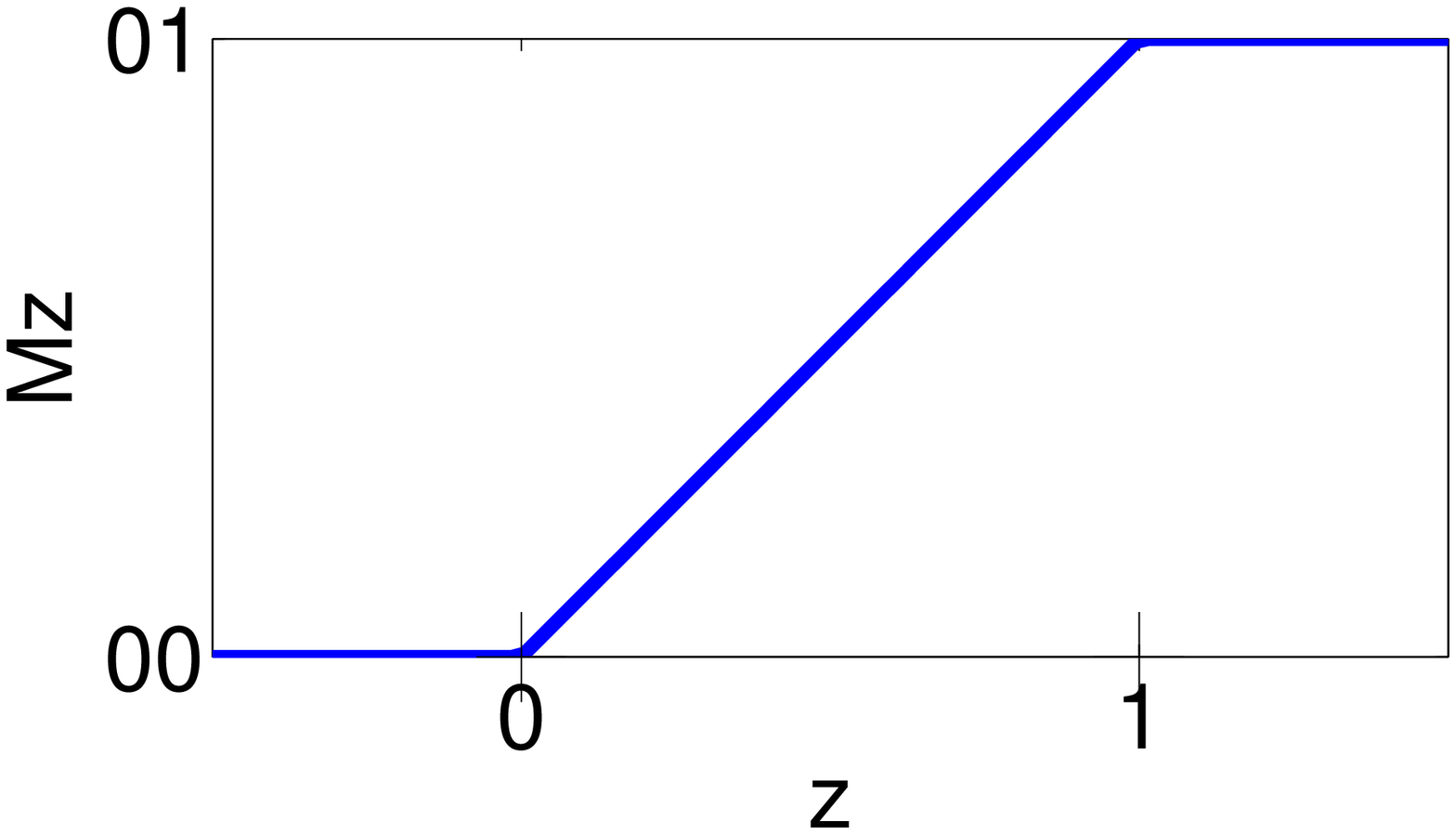}
  \end{tabular}
  \caption{\emph{Left}: three possible cases for the location of the solution for the scalar box-constrained QP in the \z\ step. \emph{Right}: the solution $M(l,u,z)$ is the median of $l$, $u$ and $z$.}
  \label{f:proxbqpZ}
\end{figure}

In order to apply ADMM to our QP~\eqref{e:proxbqp}, we make the identification $f(\x) = \frac{1}{2} \x^T \A \x - \b^T \x$, defined over the entire $\bbR^D$, and $g(\z) = \frac{\mu}{2} \smash{\norm{\z - \vv}^2}$, defined for $\bl \le \z \le \uu$, and want to minimize $f(\x) + g(\z)$ s.t.\ $\x = \z$. The augmented Lagrangian is
\begin{equation*}
  \calL(\x,\z,\blambda;\rho) = \frac{1}{2} \x^T \A \x - \b^T \x + \frac{\mu}{2} \norm{\z - \vv}^2 + \blambda^T (\x-\z) + \frac{\rho}{2} \norm{\x-\z}^2 \text{ s.t. } \bl \le \z \le \uu
\end{equation*}
and the ADMM iteration is:
\begin{align*}
  \x &\leftarrow \arg\min_{\x}{\left( \frac{1}{2} \x^T \A \x - \b^T \x + \frac{\rho}{2} \norm{\x - \z + \bzeta}^2 \right)} \\
  \z &\leftarrow \arg\min_{\z}{\left( \frac{\mu}{2} \norm{\z - \vv}^2 + \frac{\rho}{2} \norm{\x - \z + \bzeta}^2 \right) \text{ s.t. } \bl \le \z \le \uu} \\
  \bzeta &\leftarrow \bzeta + \x - \z.
\end{align*}
We now solve the steps over \x\ and \z. The step over \x\ is an unconstrained strongly convex quadratic function whose unique minimizer is given by a linear system. The step over \z\ separates over each component $z_1,\dots,z_D$ of \z\ because both the objective and the constraints are separable. For each component $z_d$, $d=1,\dots,D$, we have to minimize an upwards parabola defined in $[l_d,u_d]$, whose solution is the median $M(l_d,u_d,z^*_d)$ of $l_d$, $u_d$ and the parabola vertex $z^*_d = \smash{\frac{\mu v_d + \rho(x_d+\zeta_d)}{\mu+\rho}}$ (see fig.~\ref{f:proxbqpZ}). Finally, we obtain the following updates, which are iterated in order until convergence:
\begin{subequations}
  \label{e:admm-proxqp}
  \begin{align}
    \label{e:admm-proxqp1}
    \x &\leftarrow (\A + \rho\I)^{-1} (\b + \rho(\z-\bzeta)) \\
    \label{e:admm-proxqp2}
    \z &\leftarrow M\left( \bl,\uu,\frac{\mu \vv + \rho(\x+\bzeta)}{\mu+\rho} \right) \\
    \label{e:admm-proxqp3}
    \bzeta &\leftarrow \bzeta + \x - \z
  \end{align}
\end{subequations}
where the median applies elementwise to its vector arguments, \x\ are the primal variables, \z\ the auxiliary (consensus) variables, and \bzeta\ the scaled Lagrange multiplier estimates for \z. Moving the inequalities into the \z\ subproblem results in a very simple step simply involving elementwise thresholding operations.

The iteration~\eqref{e:admm-proxqp} has a number of advantages. It is very simple to implement. It requires no line searches and has only one user parameter, the penalty parameter $\rho$, for which a good default value can be computed (see below). The algorithm converges for any positive value of $\rho$. The steps are fast and we can stop at any time with a feasible iterate. The algorithm is particularly convenient if we have to solve $N$ proximal operators of the form~\eqref{e:proxbqp} that have the same matrix \A, since the Cholesky factor is computed just once and used by all operators. Finally, if the QP is the relaxation of a binary problem (i.e., \x\ was originally in $\{0,1\}^D$ but was relaxed to $[0,1]^D$), we need not converge to high accuracy because the final result will be binarized a posteriori.

\paragraph{Computational complexity}
\label{s:cost}

Each step in~\eqref{e:admm-proxqp} is $\calO(D)$ except for the \x\ step, which is the most costly because of the linear system. This may be solved efficiently in one of the two following ways:
\begin{itemize}
\item Preferably, by caching the Cholesky factorization of $\A + \rho\I$ (i.e., $\RR\RR^T = \A + \rho\I$ where \RR\ is upper triangular). If \A\ is dense, computing the Cholesky factor is a one-time cost of $\calO(\frac{1}{3} D^3)$, and solving the linear system is $\calO(D^2)$ by solving two triangular systems. If \A\ is large and (sufficiently) sparse, computing the Cholesky factor and solving the linear system are both possibly $\calO(D)$. One should use a good permutation to reduce fill-in. This is feasible with relatively small dimensions $D$ if using a dense \A, or with large dimensions as long as \A\ is sufficiently sparse that the Cholesky factorization does not add so much fill that it can be stored.
\item By using an iterative linear solver such as conjugate gradients \citep{NocedalWright06a}, initialized with a warm start, preconditioned, and exiting it before convergence, so as to carry out faster, inexact \x-updates. This is better for large problems where the Cholesky factorization adds too much fill.
\end{itemize}
Thus, each iteration of the algorithm is cheap. In practice, for good values of $\rho$, the algorithm quickly approaches the solution in the first iterations and then converges slowly, as is known with ADMM algorithms in general. However, since each iteration is so cheap, we can run a large number of them if high accuracy is needed. As a sample runtime, for a collection of $N=60\,000$ problems in $D=32$ variables having the same, dense matrix \A, solving all $N$ problems to high accuracy takes around $1$~s in a PC.

\paragraph{Initialization}

If the QP problem~\eqref{e:proxbqp} is itself a subproblem in a larger problem, one should warm-start the iteration of eq.~\eqref{e:admm-proxqp} from the values of \z\ and \bzeta\ in the previous outer-loop iteration. Otherwise, we can simply initialize $\z = \vv$ and $\bzeta = \0$.

\paragraph{Stopping criterion}

We stop when the relative change in \z\ is less than a set tolerance (e.g.\ $10^{-5}$). We return as approximate solution the value of \z, which is feasible by construction, while \x\ need not be (upon convergence, $\x = \z$).

\paragraph{Optimal penalty parameter $\rho$}

The speed at which ADMM converges depends significantly on the quadratic penalty parameter $\rho$ \citep{Boyd_11a}. Little work exists on how to select $\rho$ so as to achieve fastest convergence. Recently, for QPs, \citet{Ghadim_13a} suggest to use $\rho^* = \sqrt{\sigma_{\min} \sigma_{\max}}$ where $\sigma_{\min}$ and $ \sigma_{\max}$ are the smallest (nonzero) and largest eigenvalue of the matrix \A.

\paragraph{Matlab code}
\label{s:matlab}

A Matlab implementation of the algorithm is available from the author. The following Matlab code implements a basic form of the algorithm, using the Cholesky factor in the \x-update linear system, and assuming the bounds $[l,u]$ are the same for each of the $N$ quadratic programs. It applies the algorithm synchronously to all the QPs, which means each QP runs the same number of iterations. This is inefficient, because the runtime is driven by the slowest variable to converge, but the code is better vectorized in Matlab. An implementation in C should handle each QP separately, particularly in a parallel or distributed implementation.

\begin{verbatim}
function [Z,Y,X] = proxbqp(V,m,A,B,l,u,Z,Y,r,maxit,tol)

[D,N] = size(V); R = chol(A+r*eye(D,D)); Rt = R'; Zold = zeros(D,N);
for i=1:maxit
  X = R \ (Rt \ (B + r*(Z-Y)));
  Z = (m*V+r*(X+Y))/(m+r); Z(Z<l) = l; Z(Z>u) = u;
  Y = Y + X - Z;
  if max(abs(Z(:)-Zold(:))) < tol break; end; Zold = Z;
end
\end{verbatim}


\begin{thebibliography}{7}
\providecommand{\natexlab}[1]{#1}
\providecommand{\url}[1]{\texttt{#1}}
\expandafter\ifx\csname urlstyle\endcsname\relax
  \providecommand{\doi}[1]{doi: #1}\else
  \providecommand{\doi}{doi: \begingroup \urlstyle{rm}\Url}\fi

\bibitem[Boyd et~al.(2011)Boyd, Parikh, Chu, Peleato, and Eckstein]{Boyd_11a}
S.~Boyd, N.~Parikh, E.~Chu, B.~Peleato, and J.~Eckstein.
\newblock Distributed optimization and statistical learning via the alternating
  direction method of multipliers.
\newblock \emph{Foundations and Trends in Machine Learning}, 3\penalty0
  (1):\penalty0 1--122, 2011.

\bibitem[Carreira-Perpi{\~n}{\'a}n and
  Raziperchikolaei(2015)]{CarreirRaziper15a}
M.~{\'A}. Carreira-Perpi{\~n}{\'a}n and R.~Raziperchikolaei.
\newblock Hashing with binary autoencoders.
\newblock Unpublished manuscript, Jan. 2015.

\bibitem[Combettes and Pesquet(2011)]{CombetPesquet11a}
P.~L. Combettes and J.-C. Pesquet.
\newblock Proximal splitting methods in signal processing.
\newblock In H.~H. Bauschke, R.~S. Burachik, P.~L. Combettes, V.~Elser, D.~R.
  Luke, and H.~Wolkowicz, editors, \emph{Fixed-Point Algorithms for Inverse
  Problems in Science and Engineering}, Springer Series in Optimization and Its
  Applications, pages 185--212. Springer-Verlag, 2011.

\bibitem[Ghadimi et~al.(2013)Ghadimi, Teixeira, Shames, and
  Johansson]{Ghadim_13a}
E.~Ghadimi, A.~Teixeira, I.~Shames, and M.~Johansson.
\newblock Optimal parameter selection for the alternating direction method of
  multipliers {(ADMM)}: Quadratic problems.
\newblock arXiv:1306.2454 [math.OC], July~3 2013.

\bibitem[Moreau(1962)]{Moreau62a}
J.-J. Moreau.
\newblock Fonctions convexes duales et points proximaux dans un espace
  hilbertien.
\newblock \emph{C. R. Acad. Sci. Paris S{\'e}r. A Math.}, 255:\penalty0
  2897--2899, 1962.

\bibitem[Nocedal and Wright(2006)]{NocedalWright06a}
J.~Nocedal and S.~J. Wright.
\newblock \emph{Numerical Optimization}.
\newblock Springer Series in Operations Research and Financial Engineering.
  Springer-Verlag, New York, second edition, 2006.

\bibitem[Rockafellar(1976)]{Rockaf76b}
R.~T. Rockafellar.
\newblock Monotone operators and the proximal point algorithm.
\newblock \emph{SIAM J. Control and Optim.}, 14\penalty0 (5):\penalty0
  877--898, 1976.

\end{thebibliography}

\end{document}